\documentclass[11pt]{article}

\usepackage{amssymb,amsmath,amsthm}
\usepackage{mathrsfs}
\usepackage{mathtools}
\usepackage{graphicx}

\numberwithin{equation}{section}
\newtheorem{thm}{Theorem}[section]

\theoremstyle{definition}
\newtheorem{rem}[thm]{Remark}
\newtheorem{rems}[thm]{Remarks}
\let\oldproofname=\proofname
\renewcommand{\proofname}{\rm\bf{\oldproofname}}

\newcommand{\N}{\mathbb{N}}
\newcommand{\Z}{\mathbb{Z}}

\newcommand{\R}{\mathbb{R}}
\newcommand{\C}{\mathbb{C}}

\newcommand{\cA}{\mathcal{A}}
\newcommand{\cB}{\mathcal{B}}
\newcommand{\cC}{\mathcal{C}}
\newcommand{\cD}{\mathcal{D}}
\newcommand{\cE}{\mathcal{E}}

\newcommand{\cO}{\mathcal{O}}

\renewcommand{\Re}{\mathop{\mathrm{Re}}}
\renewcommand{\Im}{\mathop{\mathrm{Im}}}
\newcommand{\dd}{\,{\rm d}}
\newcommand{\D}{{\rm d}}
\renewcommand{\div}{\mathop{\mathrm{div}}\nolimits}

\newcommand{\Ri}{\mathrm{Ri}}

\newcommand{\DS}{\displaystyle}

\renewcommand{\:}{\thinspace :}
\newcommand{\?}{\thinspace ?}
\newcommand{\ggamma}{\gamma_\star}

\newdimen\texpscorrection
\newdimen\figcenter
\def\figurewithtex #1 #2 #3 #4 #5\cr{\null
  {\goodbreak\figcenter=\hsize\relax
  \advance\figcenter by -#4truecm
  \divide\figcenter by 2
  \begin{figure}[hbt]
  \vskip #3truecm\noindent\hskip\figcenter
  \includegraphics{#1}{\hskip\texpscorrection\input #2 }
  \vskip 0.8truecm{\baselineskip=0.8\baselineskip
  \noindent \vbox{\noindent {\footnotesize #5}}\par}
  \end{figure}}}
\def\point#1 #2 #3 {\rlap{\kern #1 truecm
\raise #2 truecm \hbox{#3}}}

\begin{document}

\title{Stability of Vortices in Ideal Fluids\:\\ 
the Legacy of Kelvin and Rayleigh}

\author{
{\bf Thierry Gallay}\\[1mm]
Institut Fourier, Universit\'e Grenoble Alpes\\
100 rue des Maths, 38610 Gi\`eres, France\\
{\tt Thierry.Gallay@univ-grenoble-alpes.fr}}

\date{}
\maketitle

\begin{abstract}
The mathematical theory of hydrodynamic stability started in the
middle of the 19th century with the study of model examples,
such as parallel flows, vortex rings, and surfaces of discontinuity. 
We focus here on the equally interesting case of columnar vortices,
which are axisymmetric stationary flows where the velocity
field only depends on the distance to the symmetry axis and has no
component in the axial direction. The stability of such flows was
first investigated by Kelvin in 1880 for some particular velocity
profiles, and the problem benefited from important contributions by
Rayleigh in 1880 and 1917. Despite further progress in the 20th 
century, notably by Howard and Gupta (1962), the only rigorous results 
so far are necessary conditions for instability under 
either two-dimensional or axisymmetric perturbations. This note is 
a non-technical introduction to a recent work in collaboration with 
D. Smets, where we prove under mild assumptions that columnar 
vortices are spectrally stable with respect to general three-dimensional 
perturbations, and that the linearized evolution group has a subexponential 
growth as $|t| \to \infty$. 
\end{abstract}

\section{Introduction to Hydrodynamic Stability Theory}
\label{sec1}

Hydrodynamic stability is the subdomain of fluid dynamics which
studies the stability and the onset of instability in fluid flows.
These fundamental questions were first addressed in the 19th century,
with pioneering contributions by G.~Stokes, H.~von Helmholtz,
W.~Thomson (Lord Kelvin), and J.~W.~Strutt (Lord
Rayleigh) on the theoretical side, and by O.~Reynolds on
the experimental side \cite{Da}. In early times the notion of
stability still lacked a precise mathematical definition, but 
its physical meaning was already perfectly understood, as 
can be seen from  the following quote by J.~C.~Maxwell \cite{Ha}, 
which dates back to 1873\:

\begin{quote}
``When the state of things is such that an infinitely small variation of
the present state will alter only by an infinitely small quantity the
state at some future time, the condition of the system, whether at
rest or in motion, is said to be stable; but when an infinitely small
variation in the present state may bring about a finite difference in
the state of the system in a finite time, the system is said to be
unstable.''
\end{quote}

What is exactly meant by ``infinitely small'' in this definition is
rigorously specified, for instance, in the subsequent memoir by
A.~M.~Lyapunov \cite{Ly}, which was published in 1892. The relevance
of stability questions in fluid mechanics cannot be overestimated. As
an example, in the idealized situation where the fluid is assumed to
be incompressible and inviscid, a plethora of explicit stationary
solutions are known which describe shear flows, vortices, or flows
past obstacles. However, depending on circumstances, these solutions
may or may not be observed in real life, where experimental
uncertainties, viscosity effects, and boundary conditions play an
important role. To determine the relevance of a given flow, the
stability analysis is certainly the first step to perform, but even in
an idealized framework this often leads to difficult mathematical
problems, a complete solution of which was largely out of reach in the
19th century and is still a serious challenge today.

To make the previous considerations more concrete, we analyze
in this introduction three relatively simple cases, of increasing
complexity, where stability can be discussed using the techniques
introduced by Rayleigh \cite{Ra1}. These examples are classical and
thoroughly studied in many textbooks \cite{Cha,Chu,DR,Lin,SH}, 
as well as in the excellent review article \cite{DH}. The results obtained 
for these model problems will serve as a guideline for the stability 
analysis of columnar vortices, which will be presented in 
Sections~\ref{sec2} and \ref{sec3}.

\subsection{The Rayleigh-Taylor Instability}\label{sec11}

We consider the motion of an incompressible and inviscid fluid in the
infinite strip $D = \R \times [0,L]$ with coordinates $(x,z)$, where
$x \in \R$ is the horizontal variable and $z \in [0,L]$ the vertical
variable. The state of the fluid at time $t \in \R$ is defined by the
density distribution $\rho(x,z,t) > 0$, the velocity field
$u(x,z,t) \in \R^2$, and the pressure $p(x,z,t) \in \R$. The evolution 
is determined by the density-dependent incompressible Euler equations
\begin{equation}\label{RT1}
  \partial_t \rho + u \cdot \nabla \rho = 0\,, \quad 
  \rho\bigl(\partial_t u + (u\cdot\nabla)u\bigr) \,=\, 
  - \nabla p - \rho g e_z\,, \quad \div u \,=\, 0\,,
\end{equation}
where $g$ denotes the acceleration due to gravity and $e_z$ is 
the unit vector in the (upward) vertical direction. Setting 
$u = (u_x,u_z)$, we impose the impermeability condition $u_z(x,z,t) 
= 0$ at the bottom and the top of the domain $D$, namely for $z = 0$ 
and $z = L$. 

The PDE system \eqref{RT1} has a family of stationary solutions of
the form $\rho = \bar \rho(z)$, $u = 0$, $p = \bar p (z)$, where the
density $\bar \rho$ is an arbitrary function of the vertical
coordinate $z$, and the associated pressure is determined (up to an
irrelevant additive constant) by the hydrostatic balance
$\bar p'(z) = -\bar \rho(z)g$. To study the stability of the
equilibrium $(\bar \rho,0,\bar p)$, we consider perturbed solutions
of the form
\[
  \rho(x,z,t) \,=\, \bar \rho(z) + \tilde \rho(x,z,t)\,, \quad
  u(x,z,t) \,=\, \tilde u(x,z,t)\,, \quad
  p(x,z,t) \,=\, \bar p(z) + \tilde p(x,z,t)\,. 
\]
Inserting this Ansatz into \eqref{RT1} and neglecting all 
quadratic terms in $(\tilde \rho,\tilde u)$, we obtain the 
{\em linearized} equations for the perturbations $(\tilde \rho,
\tilde u,\tilde p)$\:
\begin{equation}\label{RT2}
  \begin{array}{l}
   \bar \rho(z) \partial_t \tilde u_x \,=\, -\partial_x \tilde p\,,\\[1mm]   
   \bar \rho(z) \partial_t \tilde u_z \,=\, -\partial_z \tilde p 
   -\tilde \rho g\,, 
  \end{array}\qquad\quad
  \begin{array}{r}
   \partial_t \tilde \rho + \bar\rho'(z)\tilde u_z \,=\, 0\,, \\[1mm]
   \partial_x \tilde u_x + \partial_z \tilde u_z \,=\, 0\,.
  \end{array}   
\end{equation}

\begin{rem}\label{rem:lin}
It is not obvious at all that considering the linearized perturbation 
equations \eqref{RT2} is sufficient, or even appropriate, to determine 
the stability of stationary solutions to \eqref{RT1}. In fact the 
validity of Lyapunov's linearization method in the context of fluid 
mechanics is a difficult question \cite{Yu}, which is the object of 
ongoing research. In particular, for ideal fluids, there is no general 
result asserting that a linearly stable equilibrium is actually stable 
in the sense of Lyapunov. However, if the linearized system is 
exponentially unstable, for instance due to the existence of an 
eigenvalue with nonzero real part, it is often possible 
to conclude that the equilibrium under consideration is unstable, 
see \cite{BGS,FSV,Li2,VF} for a few results in this direction. 
To summarize, the linearization approach may be useful to detect
exponential instabilities, but stability results have to be established
by a different approach, for instance (in two space dimensions) using
variational techniques \cite{Arn1,Arn2}
\end{rem}

The linearized equations \eqref{RT2} are invariant under
translations in the horizontal direction, so that we can use a Fourier
transform to reduce the number of independent variables. A further
simplification is made by restricting our attention to {\em
eigenfunctions} of the linearized operator. In other words, we
consider solutions of \eqref{RT2} of the particular form
\begin{equation}\label{RT3}
  \tilde \rho(x,z,t) \,=\, \rho(z)\,e^{ikx}\,e^{st}\,, \quad
  \tilde u(x,z,t) \,=\, u(z)\,e^{ikx}\,e^{st}\,, \quad
  \tilde p(x,z,t) \,=\, p(z)\,e^{ikx}\,e^{st}\,,
\end{equation}
where $k \in \R$ is the horizontal wave number and $s \in \C$ is 
the spectral parameter. The representation \eqref{RT3} transforms 
the linearized equations \eqref{RT2} into an ODE system\:
\begin{equation}\label{RT4}
  \begin{array}{l}
   \bar \rho(z) s  u_x \,=\, -ik  p\,,\\[1mm]   
   \bar \rho(z) s  u_z \,=\, -\partial_z  p 
   - \rho g\,, 
  \end{array}\qquad\quad
  \begin{array}{r}
   s  \rho + \bar\rho'(z) u_z \,=\, 0\,, \\[1mm]
   ik  u_x + \partial_z  u_z \,=\, 0\,,
  \end{array}   
\end{equation}
which (if $s \neq 0$) can in turn be reduced to a single equation 
for the vertical velocity profile $u_z$\:
\begin{equation}\label{RT5}
   -\partial_z \bigl(\bar \rho(z)\partial_z u_z\bigr) + 
  k^2 \bar \rho(z) u_z - \frac{k^2 g}{s^2}\, \bar\rho'(z)u_z 
  \,=\, 0\,, \qquad z \in [0,L]\,.
\end{equation}
By construction, the values of the spectral parameter
$s \in \C\setminus \{0\}$ for which the ODE \eqref{RT5} has a
nontrivial solution $u_z$ satisfying the boundary conditions
$u_z(0) = u_z(L) = 0$ are {\em eigenvalues} of the linearized operator
\eqref{RT2} in the Fourier subspace indexed by the horizontal
wavenumber $k \in \R$. Spectral stability is obtained if all
eigenvalues are purely imaginary, whereas the existence of an
eigenvalue $s \in \C$ with $\Re(s) \neq 0$ implies exponential
instability of the linearized system in positive or negative times.

\begin{rems}\label{rem:spec}\quad\\
{\bf 1.} The Fourier transform reduces the linearized equations to 
a one-dimensional PDE system in the bounded domain $[0,L]$, but 
this does not immediately imply that the spectrum of the full 
linearized operator is the union of the point spectra obtained 
for all values of the horizontal wavenumber $k \in \R$. So, 
even if one can prove that the eigenvalues are purely imaginary
for all $k \in \R$, an additional argument is needed to verify 
that the full linearized operator has indeed no spectrum 
outside the imaginary axis. This rather technical issue will not 
be discussed further in this introduction, but we shall come 
back to it in Section~\ref{sec3}. \\[1mm]
{\bf 2.} In the literature, the Rayleigh-Taylor equation \eqref{RT5} is 
often derived in the Boussinesq approximation, which consists 
in neglecting the variations of the density profile $\bar \rho(z)$ 
everywhere except in the buoyancy term. This gives the simplified
eigenvalue equation
\begin{equation}\label{RT6}
  -\partial_z^2 u_z + k^2 \Bigl(1 + \frac{N(z)^2}{s^2}\Bigr)u_z  
   \,=\, 0\,, \qquad \hbox{where}\quad 
   N(z)^2 \,=\, - \frac{g \bar\rho'(z)}{\rho(z)}\,.
\end{equation}
When $\bar \rho'(z) < 0$, the real number $N(z)$ is called the 
Brunt-V\"ais\"al\"a frequency. This is the (maximal) oscillation
frequency of gravity waves inside a stably stratified fluid. 
\end{rems}

Assume that, for some $k \in \R$ and some $s \in \C \setminus\{0\}$, 
the ODE \eqref{RT5} has a nontrivial solution $u_z$ satisfying
the boundary conditions $u_z(0) = u_z(L) = 0$. Multiplying 
both sides of \eqref{RT5} by the complex conjugate $\bar u_z$ 
and integrating over the vertical domain $[0,L]$, we obtain
the integral identity
\begin{equation}\label{RT7}
  \int_0^L \bar \rho(z) |\partial_z u_z|^2 \dd z + 
  k^2  \int_0^L \bar \rho(z) |u_z|^2 \dd z - 
  \frac{k^2 g}{s^2} \int_0^L \bar \rho'(z) |u_z|^2 \dd z 
  \,=\, 0\,.
\end{equation}
The first two terms in \eqref{RT7} being real and positive, equality
can hold only if the third term is real and negative. Thus we must
have $k \neq 0$ and $\Im(s^2) = 0$, namely $s \in \R$ or $s \in
i\R$. Now, if we assume that the fluid is {\em stably stratified}, in 
the sense that $\bar \rho'(z) \le 0$ for all $z \in [0,L]$, the 
last term in \eqref{RT7} is positive only if $s^2 < 0$, which means
that $s \in i\R$. Under this assumption, we conclude that all 
eigenfunctions of the form \eqref{RT3} with $k \in \R$ correspond
to eigenvalues $s$ on the imaginary axis, so that the equilibrium 
$(\bar \rho,0,\bar p)$ of \eqref{RT1} is {\em spectrally stable}, 
up to the technical issue mentioned in Remark~1.2.1. 

On the other hand, if $\bar \rho'(z) > 0$ for some $z \in [0,L]$, a
nice argument due to Synge \cite{Sy} shows that, for any $k \neq 0$,
the Rayleigh equation has a nontrivial solution $u_z$ (satisfying the
boundary conditions) for a sequence of real eigenvalues $s_n \to 0$.
The equilibrium $(\bar \rho,0,\bar p)$ of \eqref{RT1} is thus 
spectrally unstable. Summarizing, the stability of the rest state 
$u = 0$ in stratified ideal fluids is reasonably understood, 
in the sense that Rayleigh's approach provides a {\em necessary and 
sufficient} condition for spectral stability in that case. 

\subsection{Shear Flows in Homogeneous Fluids}\label{sec12}

For the same equations \eqref{RT1} in the domain $D$, we now consider
a different family of equilibria, namely shear flows of the form
$\rho = 1$, $u = U(z)e_x$, $p = 0$, where the horizontal velocity
profile $U$ is an arbitrary function. For the moment, we assume that
the fluid is homogeneous and only allow for perturbations of the
velocity field.  The perturbed solutions thus take the form
\[
  \rho(x,z,t) \,=\, 1\,, \quad u(x,z,t) \,=\, U(z)e_x + \tilde u(x,z,t)\,, 
  \quad p(x,z,t) \,=\, \tilde p(x,z,t)\,, 
\]
and the linearized equations become
\begin{equation}\label{SF1}
  \begin{array}{r}
  \partial_t \tilde u_x + U(z)\partial_x \tilde u_x + U'(z)
  \tilde u_z \,=\, -\partial_x \tilde p\,, \\[1mm]
  \partial_t \tilde u_z + U(z)\partial_x \tilde u_z \,=\, 
  -\partial_z \tilde p\,,
  \end{array} \qquad\quad
  \partial_x \tilde u_x + \partial_z \tilde u_z \,=\, 0\,.
\end{equation}
As before, we suppose that $\tilde u(x,z,t) = u(z)\,e^{ikx}\,e^{st}$ 
and $\tilde p(x,z,t) = p(z)\,e^{ikx}\,e^{st}$ for some $k \in \R$ and 
some $s \in \C$. The functions $u,p$ are solutions of the 
ODE system
\begin{equation}\label{SF2}
  \gamma(z) u_x + U'(z) u_z \,=\, -ikp\,, \quad
  \gamma(z) u_z \,=\, -\partial_z p\,, \quad
  iku_x + \partial_z u_z \,=\, 0\,,
\end{equation}
where $\gamma(z) = s + ikU(z)$ is the symbol of the material 
derivative $\partial_t + U(z)\partial_x$. This function plays 
an important role in the stability analysis, as it incorporates 
the spectral parameter $s$. 

Since we are interested in detecting potential 
instabilities, we assume in what follows that $\Re(s) \neq 0$, 
which implies in particular that $\gamma(z) \neq 0$ for all 
$z \in [0,L]$. Under this hypothesis, we can reduce the ODE 
system \eqref{SF2} to the following scalar equation for the 
vertical velocity\:
\begin{equation}\label{SF3}
   -\partial_z^2 u_z + k^2 u_z + \frac{ikU''(z)}{\gamma(z)}\,u_z 
  \,=\, 0\,, \qquad z \in [0,L]\,.
\end{equation}
This equation looks simpler than \eqref{RT5}, but is in fact 
substantially harder to analyze. If $u_z$ is a nontrivial solution
satisfying the boundary conditions, we have Rayleigh's 
identity
\begin{equation}\label{SF4}
  \int_0^L |\partial_z u_z|^2 \dd z + k^2  \int_0^L |u_z|^2 \dd z + 
  ik \int_0^L \frac{U''(z)}{\gamma(z)}\,|u_z|^2 \dd z \,=\, 0\,,
\end{equation}
which can be satisfied only if $k \neq 0$ and if $U''(z)$ is 
not identically zero. Under these assumptions, the 
imaginary part of \eqref{SF4} gives the useful relation
\begin{equation}\label{SF5}
  \Re(s) \int_0^L \frac{U''(z)}{|\gamma(z)|^2}\,|u_z|^2 \dd z \,=\, 0\,.
\end{equation}
If $U''(z)$ does not change sign on $[0,L]$, the integral in 
\eqref{SF5} is nonzero, which contradicts our assumption that
$\Re(s) \neq 0$. This gives Rayleigh's {\em inflection point 
criterion} \cite{Ra1}\: a necessary condition for the shear flow 
with velocity profile $U(z)$ to be (spectrally) unstable is that 
the function $z \mapsto U''(z)$ changes sign on the interval 
$[0,L]$. 

Rayleigh's inflection point criterion is not sharp, and can 
be improved somehow by exploiting both the real and imaginary 
parts of identity \eqref{SF4}, see \cite{Fj}. However, 
surprisingly enough, it seems difficult to formulate a necessary 
and sufficient stability condition for shear flows, even in 
the ideal case considered here. An instructive example is 
Kolmogorov's flow $U(z) = \sin(z - L/2)$, which is known to 
be stable if and only if $L \le \pi$ \cite{DH,Li1}, although both 
Rayleigh's and Fj\o rtoft's criteria allow for a possible 
instability for any $L > 0$. In fact, the origin of inertial instabilities 
in shear flows seems only partially understood from a physical point 
of view, see \cite{BM}. 

\subsection{Shear Flows in Stratified Fluids}\label{sec13}

Following the same approach as in the previous paragraphs, we now
analyze the stability of shear flows in (stably) stratified fluids. 
We consider the Euler equations \eqref{RT1} in the vicinity of a
stationary solution of the form $\rho = \bar \rho(z)$, $u = U(z) e_x$,
$p = \bar p(z)$, where $\bar p'(z) = -\bar \rho(z)g$ (hydrostatic
balance).  The perturbed solutions are written in the form
\[
  \rho(x,z,t) \,=\, \bar \rho(z) + \tilde \rho(x,z,t)\,, \quad
  u(x,z,t) \,=\, U(z) e_x + \tilde u(x,z,t)\,, \quad
  p(x,z,t) \,=\, \bar p(z) + \tilde p(x,z,t)\,,
\]
so that the linearized equations become
\begin{equation}\label{TG1}
  \begin{array}{l}
  \bar \rho(z)\bigl(\partial_t \tilde u_x + U(z)\partial_x 
  \tilde u_x + U'(z) \tilde u_z\bigr) \,=\, -\partial_x \tilde p\,,\\[1mm]   
  \bar \rho(z)\bigl(\partial_t \tilde u_z + U(z)\partial_x 
  \tilde u_z\bigr) \,=\, -\partial_z \tilde p -\tilde \rho g\,, 
  \end{array}\qquad
  \begin{array}{r}
  \partial_t \tilde \rho + U(z)\partial_x \tilde\rho + \bar\rho'(z)
  \tilde u_z \,=\, 0\,, \\[1mm]
   \partial_x \tilde u_x + \partial_z \tilde u_z \,=\, 0\,.
  \end{array}   
\end{equation}
For perturbations of the form \eqref{RT3}, we arrive at the ODE 
system 
\begin{equation}\label{TG2}
  \begin{array}{l}
   \bar \rho(z)\bigl(\gamma(z)u_x + U'(z)u_z\bigr) \,=\, -ik  p\,,\\[1mm]   
   \bar \rho(z) \gamma(z) u_z \,=\, -\partial_z  p - \rho g\,, 
  \end{array}\qquad\quad
  \begin{array}{r}
   \gamma(z) \rho + \bar\rho'(z) u_z \,=\, 0\,, \\[1mm]
   ik  u_x + \partial_z  u_z \,=\, 0\,,
\end{array}   
\end{equation}
where $\gamma(z) = s + ikU(z)$ is the spectral function. If we assume 
that $\Re(s) \neq 0$, so that $\gamma(z) \neq 0$, we can reduce 
the system \eqref{TG2} to the {\em Taylor-Goldstein equation}
\begin{equation}\label{TG3}
  -\partial_z \bigl(\bar \rho(z)\partial_z u_z\bigr) + 
  k^2 \bar \rho(z) u_z + \frac{ik}{\gamma(z)}\bigl(\bar\rho U'
  \bigr)'(z) u_z - \frac{k^2 g}{\gamma(z)^2}\,\bar\rho'(z)u_z 
  \,=\, 0\,, \qquad z \in [0,L]\,.
\end{equation}
Note that we recover the Rayleigh-Taylor equation \eqref{RT5}
by setting $U = 0$, hence $\gamma(z) = s$, in \eqref{TG2}. 
Similarly, \eqref{TG2} reduces to the Rayleigh stability equation 
\eqref{SF3} when $\bar\rho = 1$. 

The original approach of Rayleigh does not give much information 
on the solutions of \eqref{TG3}. If $u_z$ is a nontrivial solution 
satisfying the boundary conditions, it is difficult to exploit 
the integral identity
\begin{equation}\label{TG4}
  \int_0^L \biggl\{\bar \rho(z) |\partial_z u_z|^2 + 
  k^2 \bar \rho(z) |u_z|^2 + ik \frac{(\bar\rho U')'(z)}{\gamma(z)}\,|u_z|^2
  - \frac{k^2 g}{\gamma(z)^2} \bar \rho'(z) |u_z|^2\biggr\} \dd z 
  \,=\, 0\,,
\end{equation}
because the real or imaginary parts of the last two terms in the integrand 
have no obvious sign. A solution to this problem was found by Miles 
\cite{Mil} and Howard \cite{How} in the early 60's. Following the 
elegant approach of \cite{How}, we perform the change of variables
\[
  u_z(z) \,=\, \gamma(z)^{1/2} v_z(z)\,, \qquad \hbox{where}\quad 
  \gamma(z) \,=\, s + ikU(z)\,.
\]
The new function $v_z$ satisfies the modified ODE 
\begin{equation}\label{TG5}
  -\partial_z \bigl(\bar \rho(z)\gamma(z)\partial_z v_z\bigr) +
  k^2 \bar \rho(z)\gamma(z) v_z + \frac{ik}{2}\,(\bar \rho U')'(z)
  \,v_z + \Bigl(\frac{\bar \rho \gamma'^2}{4\gamma} - \frac{k^2 g 
  \bar\rho'}{\gamma}\Bigr)(z)\,v_z \,=\, 0\,.
\end{equation}
If $v_z$ is a nontrivial solution satisfying the boundary conditions
$v_z(0) = v_z(L) = 0$, we multiply both sides of \eqref{TG5} by 
the complex conjugate $\bar v_z$ and integrate over the domain $[0,L]$. 
After taking the real part, we obtain the useful identity 
\begin{equation}\label{TG6}
   \Re(s)\int_0^L \biggl\{\bar \rho(z)\bigl(|\partial_z v_z|^2 + k^2 
  |v_z|^2\bigr) + \frac{k^2 \bar\rho(z) U'(z)^2}{|\gamma(z)|^2}
  \Bigl(\Ri(z) - \frac14\Bigr)\,|v_z|^2 \biggr\}\dd z \,=\, 0\,,
\end{equation}
where $\Ri(z)$ is the (local) {\em Richardson number} defined 
by
\[
  \Ri(z) \,=\, \frac{-\bar\rho'(z)\,g}{\bar\rho(z)}\,\frac{1}{U'(z)^2}
  \,=\, \left(\frac{N(z)}{U'(z)}\right)^2\,.
\]
We assume here that $\bar\rho'(z) \le 0$ (stable stratification), 
so that $\Ri(z) \ge 0$, and we denote by $N(z)$ the Brunt-V\"ais\"al\"a 
frequency \eqref{RT6}. 

The Richardson number compares the stabilizing effect of the
stratification, measured by the frequency $N$ of the gravity waves, to
the potentially destabilizing effect of the shear flow, which may be
proportional to the velocity gradient $U'$ \cite{DH}.  Clearly,
equality \eqref{TG6} cannot hold if $\Ri(z) \ge 1/4$ for all
$z \in [0,L]$, because the integrand is then positive while we assumed
that $\Re(s) \neq 0$. This gives the celebrated {\em Miles-Howard
  criterion}\: a shear flow in a stratified fluid is spectrally stable
if the Richardson number is greater than or equal to $1/4$ everywhere
in the fluid. The threshold value $1/4$ is known to be sharp, in the
sense that it cannot be replaced by any smaller real number. However,
the Miles-Howard criterion itself is by no means sharp\: if
$\bar \rho = 1$, any shear flow without inflection point is spectrally
stable by Rayleigh's criterion, although $\Ri(z) \equiv 0$ in that
case.

\begin{rem}\label{rem:squire}
So far we concentrated on the two-dimensional case, but it is 
also instructive to investigate the stability of shear flows 
with respect to three-dimensional perturbations. In that case, 
we work in the domain $D' = \R^2 \times [0,L]$ with coordinates 
$(x,y,z)$, and consider perturbations that are plane waves with 
horizontal wave vector $k = (k_1,k_2) \in \R^2$. For instance, 
the three-dimensional velocity field takes the form
\[
  u(x,y,z,t) \,=\, U(z) e_x + u(z)\,e^{i(k_1 x+k_2  y)}\,e^{\sigma t}\,, 
\]
where $\sigma \in \C$ is the spectral parameter. Using a similar
Ansatz for the density and the pressure, it is easy to derive the 3D
perturbation equations which generalize \eqref{TG2}. Now, in the
homogeneous case where $\rho \equiv 1$, a well-know result due to
Squire \cite{Sq} shows that, if the 3D perturbation equations have a
nontrivial solution for some $k_1,k_2 \neq 0$ and some $\sigma \in \C$
with $\Re(\sigma) \neq 0$, then the 2D perturbation equations
\eqref{SF2} also have a nontrivial solution with
$k = (k_1^2 + k_2^2)^{1/2}$ and $s = (k/k_1)\sigma$.  Note that
$|\Re(s)| > |\Re(\sigma)|$, so that the most unstable modes are always
two-dimensional; in other words, it is sufficient to consider the 2D
case to detect potential instabilities. A similar result holds in the
general situation where the fluid is stratified \cite{DH}, but in that
case Squire's transformation also affects the acceleration due to
gravity, replacing $g$ by the larger quantity $(k^2/k_1^2)g$. This
means that, to any unstable 3D mode, there corresponds a more unstable
2D mode {\em in a stronger gravitational field}. Therefore, unless the
fluid is stably stratified, this result does not imply that the most
unstable modes are necessarily two-dimensional.
\end{rem}

\section{Classical Stability Results for Vortices in Ideal Fluids}
\label{sec2}

We now discuss our main topic, namely the stability of a family of
axisymmetric stationary solutions to the three-dimensional Euler
equations which describe steady vortex columns. For symmetry reasons,
it is convenient to introduce cylindrical coordinates $(r,\theta,z)$,
and to decompose the velocity of the fluid as $u = u_r e_r + u_\theta e_\theta 
+ u_z e_z$, where $e_r$, $e_\theta$, $e_z$ are unit vectors in the radial, 
azimuthal, and vertical directions, respectively. Assuming that the fluid 
density is constant and equal to one, the Euler equations become
\begin{equation}\label{Eulercyl}
  \begin{split}
  \partial_t u_r + (u\cdot\nabla)u_r - \frac{u_\theta^2}{r} \,&=\, 
  -\partial_r p\,, \\
  \partial_t u_\theta + (u\cdot\nabla)u_\theta + \frac{u_r u_\theta}{r} \,&=\, 
  -\frac1r \partial_\theta p\,, \\
  \partial_t u_z + (u\cdot\nabla)u_z \,&=\, -\partial_z p\,,
  \end{split}
\end{equation}
where $u\cdot \nabla = u_r \partial_r  + \frac1r u_\theta \partial_\theta  
+ u_z \partial_z$. In addition, we have the incompressibility condition
\begin{equation}\label{incomp}
  \div u \,\equiv\, \frac1r\partial_r (ru_r) + \frac1r \partial_\theta u_\theta 
  + \partial_z u_z \,=\, 0\,. 
\end{equation}

Columnar vortices are stationary solutions of \eqref{Eulercyl}, 
\eqref{incomp} of the form
\begin{equation}\label{column}
  u \,=\, V(r) \,e_\theta\,,   \qquad p \,=\, P(r)\,,
\end{equation}
where $V : \R_+ \to \R$ is an arbitrary velocity profile, and the
associated pressure $P : \R_+ \to \R$ is determined, up to an
irrelevant additive constant, by the centrifugal balance
$rP'(r) = V(r)^2$. For the moment, we only assume that $V$ is a
piecewise differentiable function, and that the vortex \eqref{column}
is localized in the sense that $V(r) \to 0$ as $r \to \infty$, but
more restrictive assumptions will be added later. We introduce
the angular velocity $\Omega$ and the vorticity $W$, which are 
defined as follows\:
\begin{equation}\label{OmW}
  \Omega(r) \,=\, \frac{V(r)}{r}\,, \qquad 
  W(r) \,=\, \frac{1}{r}\,\frac{\D}{\D r}\bigl(r V(r)\bigr) \,=\, 
  r \Omega'(r) + 2 \Omega(r)\,.
\end{equation}
Without loss of generality, we normalize the vortex so that $\Omega(0) = 1$, 
hence $W(0) = 2$. Typical examples we have in mind are

\medskip
\noindent $\quad\bullet$ The {\em Rankine vortex}\: 
$\DS~\Omega(r) \,=\, \begin{cases} 1 & \hbox{if} \quad r \le 1\,, 
\\ r^{-2} & \hbox{if} \quad r \ge 1\,, \end{cases}\,$  
$\DS \quad W(r) \,=\, \begin{cases} 2 & \hbox{if} \quad r < 1\,, \\ 
0 & \hbox{if} \quad r > 1\,. \end{cases}$

\medskip
\noindent $\quad\bullet$ the {\em Lamb-Oseen vortex}\: $\DS 
~\Omega(r) \,=\, \frac{1}{r^2}\Bigl(1 - e^{-r^2}\Bigr)\,$, $\DS 
\quad W(r) \,=\, 2\,e^{-r^2}\,$. 

\medskip
\noindent $\quad\bullet$ the {\em Kaufmann-Scully vortex}\: $\DS 
~\Omega(r) \,=\, \frac{1}{1+r^2}\,$, $\DS \quad W(r) \,=\, 
\frac{2}{(1+r^2)^2}\,$.

\medskip
To study the stability of the vortex \eqref{column}, we consider 
perturbed solutions of the form
 \[
  u(r,\theta,z,t) \,=\, V(r) \,e_\theta + \tilde u(r,\theta,z,t)\,, 
  \qquad p(r,\theta,z,t) \,=\, P(r) + \tilde p(r,\theta,z,t)\,.
\]
This leads to the linearized evolution equations 
\begin{equation}\label{upert}
  \begin{split}
  \partial_t \tilde u_r + \Omega(r) \partial_\theta \tilde u_r - 2 \Omega(r) 
   \tilde u_\theta \,&=\, -\partial_r \tilde p\,, \\
  \partial_t \tilde u_\theta + \Omega(r) \partial_\theta \tilde u_\theta + 
   W(r) \tilde u_r \,&=\, -\frac1r \partial_\theta \tilde p\,, \\
  \partial_t \tilde u_z + \Omega(r) \partial_\theta  \tilde u_z \,&=\, 
  -\partial_z \tilde p\,,
  \end{split}
\end{equation}
where the pressure is determined so that the velocity perturbation 
remains divergence-free. Taking the divergence of both sides in 
\eqref{upert}, we obtain for $\tilde p$ the second order elliptic 
equation
\begin{equation}\label{pressure}
  -\partial_r^* \partial_r \tilde p -\frac{1}{r^2}\partial_\theta^2 \tilde p 
  -\partial_z^2 \tilde p
  \,=\, 2 \bigl(\partial_r^* \Omega\bigr) \partial_\theta \tilde u_r 
  - 2 \partial_r^* \bigl(\Omega\,\tilde u_\theta\bigr)\,,
\end{equation}
where we introduced the shorthand notation $\partial_r^* = \partial_r + \frac1r$. 

System \eqref{upert} was first studied by Kelvin \cite{Ke} for some
particular velocity profiles. In \cite{Ra2}, Rayleigh drew an
interesting analogy between columnar vortices and shear flows in
stratified fluids, on the basis of which he obtained a sufficient
condition for stability with respect to axisymmetric perturbations. 
Further progress was made in the 20th century, notably by Howard and
Gupta \cite{HG}, and the state of the art is reviewed in textbooks 
on vortex dynamics \cite{AKO,Sa} or hydrodynamic stability 
\cite{Cha,DR}. In this section we give a brief account of these 
classical developments, and we postpone the presentation of our 
own results to Section~\ref{sec3}.
 
\subsection{Normal Mode Analysis}\label{sec21}

As the coefficients in \eqref{upert} only depend on the distance $r$ to the 
symmetry axis, we can reduce the number of independent variables by using 
a Fourier series decomposition in the angular variable $\theta$ and a Fourier 
transform in the vertical coordinate $z$.  Moreover, as in 
Sections~\ref{sec11}--\ref{sec13}, we focus our attention to the eigenvalues 
of the linearized operator. We thus consider velocities and pressures of 
the following form
\begin{equation}\label{Four}
  \tilde u(r,\theta,z,t) \,=\, u(r)\,e^{i(m\theta+kz)}\,e^{st}\,,
  \qquad 
  p(r,\theta,z,t) \,=\, p(r)\,e^{i(m\theta+kz)}\,e^{st}\,,
\end{equation}
where $m \in \Z$ is the angular Fourier mode, $k \in \R$ is the
vertical wave number, and $s \in \C$ is the spectral parameter. 
The velocity $u = (u_r,u_\theta,u_z)$ and the pressure $p$ in \eqref{Four}
satisfy the ODE system
\begin{equation}\label{upertODE}
  \gamma(r) u_r - 2 \Omega(r) u_\theta \,=\, -\partial_r p\,, \quad
  \gamma(r) u_\theta + W(r) u_r \,=\, -\frac{im}{r} p\,, \quad
  \gamma(r) u_z \,=\, -ik p\,,
\end{equation}
where $\gamma(r) = s + im\Omega(r)$ is the spectral function. The 
incompressibility condition becomes
\begin{equation}\label{incompODE}
  \frac1r\partial_r (ru_r) + \frac{im}{r} u_\theta + ik u_z \,=\, 0\,. 
\end{equation}
If $(m,k) \neq (0,0)$ it is possible to reduce the system 
\eqref{upertODE}, \eqref{incompODE} to a scalar equation for the 
radial velocity $u_r$, by eliminating the pressure $p$ and the 
velocity components $u_\theta$, $u_z$, see \cite[Section 15]{DR} or 
\cite[Section~2]{GS1}. After straightforward calculations, we obtain 
the second order ODE
\begin{equation}\label{eigscalar}
   -\partial_r \biggl(\frac{r^2 \partial_r^* u_r}{m^2 + k^2 r^2}\biggr) 
   + \biggl\{1 + \frac{imr}{\gamma(r)}\partial_r \Bigl(\frac{W(r)}{m^2+k^2r^2}
   \Bigr) + \frac{1}{\gamma(r)^2}\frac{k^2 r^2 
  \Phi(r)}{m^2 + k^2r^2}\biggr\}u_r \,=\, 0\,,
\end{equation}
where $\partial_r^* = \partial_r + \frac1r$ and $\Phi(r) = 2\Omega(r) W(r)$ 
is the Rayleigh function. This equation is well defined if $\gamma(r) \neq 0$ 
for all $r > 0$, which is the case if $\Re(s) \neq 0$ or, more generally, 
if $s \neq -imb$ for all $b$ in the range of the angular velocity $\Omega$.
Eigenvalues of the linearized operator correspond to those values of 
the spectral parameter $s$ for which equation \eqref{eigscalar} has 
a nontrivial solution $u_r$ that is regular at the origin and decays 
to zero as $r \to \infty$. 

It is instructive to notice that the stability equation \eqref{eigscalar}
has a very similar structure as the Taylor-Goldstein equation 
\eqref{TG3}. Both are second order Schr\"odinger equations involving  
a complex-valued potential which is a polynomial of degree two in the 
inverse spectral function $1/\gamma$. The coefficient of $1/\gamma(r)^2$ 
in \eqref{eigscalar} is proportional to the Rayleigh function $\Phi$, 
and corresponds to the buoyancy term involving $-k^2 g\bar\rho'$ in 
\eqref{TG3}. Similarly, the coefficient of $1/\gamma(r)$ in \eqref{eigscalar}
is proportional to the vorticity $W$ and its derivative, and corresponds 
to the inertial term involving $ik(\bar\rho U')'$ in \eqref{TG3}. 
This analogy is grounded in deep physical reasons, which are explained 
in the pioneering work of Rayleigh \cite{Ra2}. It gives hope that 
the stability equation \eqref{eigscalar} can be analyzed using the 
techniques that were developed for shear flows, but we shall see that
additional difficulties arise in the case of columnar vortices. 

\subsection{Kelvin's Vibration Modes}\label{sec22}

When the spectral parameter $s$ is purely imaginary, the stability 
equation \eqref{eigscalar} has real-valued coefficients and can be studied
using classical methods such as Sturm-Liouville theory. If $m \neq 0$, 
it is convenient to set $s = -imb$ for some $b \in \R$, so that 
$\gamma(r) = im(\Omega(r)-b)$. In that case, the equation satisfied 
by the radial velocity $u_r$ becomes
\begin{equation}\label{eigreal}
   -\partial_r \biggl(\frac{r^2 \partial_r^* u_r}{m^2 + k^2 r^2}\biggr) 
   + \biggl\{1 + \frac{r}{\Omega(r){-}b}\partial_r \Bigl(\frac{W(r)}{m^2+k^2r^2}
   \Bigr) - \frac{1/m ^2}{(\Omega(r){-}b)^2} \frac{k^2 r^2  \Phi(r)}{m^2 
   + k^2r^2}\biggr\}u_r \,=\, 0\,.
\end{equation}
This equation is well-posed if the spectral parameter $b$ does not 
belong to the range of the angular velocity $\Omega$, so that 
$\Omega(r) - b \neq 0$ for all $r > 0$. 

In the particular case of Rankine's vortex, for which the vorticity distribution 
$W$ is piecewise constant, Kelvin \cite{Ke} observed that the stability equation 
can be explicitly solved in terms of modified Bessel functions in both regions 
$r < 1$ and $r > 1$. In the generic case where $k \neq 0$, matching 
conditions at $r = 1$ lead to the ``dispersion relation" 
\begin{equation}\label{dispersion}
  \frac{I_m'(\beta)}{\beta I_m(\beta)} + \frac{2}{(1-b)\beta^2}
  \,=\, \frac{K_m'(k)}{k K_m(k)}\,, \qquad \hbox{where}\quad
  \beta^2 \,=\, k^2\Bigl(1 - \frac{4}{m^2(1-b)^2}\Bigr)\,.
\end{equation}
Here $I_m, K_m$ are modified Bessel functions of order $m$ of the
first and second kind, respectively. Those values of $b \neq 1$ for
which \eqref{dispersion} holds give purely imaginary eigenvalues of
the linearized operator, which correspond to periodic oscillations of
the columnar vortex. A careful analysis \cite{Ke} reveals that the
relation \eqref{dispersion} is satisfied for a decreasing sequence
$b_n \to 1$, and also for an increasing sequence $b_n' \to 1$, all
solutions being contained in the interval $|b-1| \le 2/|m|$. So, for
any $m \neq 0$ and $k \neq 0$, Kelvin established the existence of an
infinite sequence of purely imaginary eigenvalues for the linearized
operator at Rankine's vortex. He was confident that the whole spectrum
could be obtained in that way \cite{Ke}\:

\begin{quote}
``All possible simple harmonic vibrations are thus found\: and
summation, after the manner of Fourier, for different values of
$[m,k]$, with different amplitudes and different epochs, gives every
possible motion, deviating infinitely little from the undisturbed
motion in circular orbits.''
\end{quote}

Unfortunately, the above claim is not substantiated by any argument in
Kelvin's paper. Nevertheless, in the case of Rankine's vortex, one can
show that the linearized operator has no eigenvalue outside the
imaginary axis, so that the whole spectrum can indeed be obtained as
demonstrated by Kelvin, see \cite[Section~6.2]{GS1}.

The situation is different for a vortex with smooth angular velocity
profile, as is the case for the Lamb-Oseen or the Kaufmann-Scully
vortex. Assuming that $\Omega(0) =1$ and $\Omega'(r) < 0$ for $r > 0$,
it can be proved that, if $m \neq 0$ and $k \neq 0$, there exists a
decreasing sequence $b_n \to 1$ of values of the spectral parameter
for which the eigenvalue equation \eqref{eigreal} has a nontrivial
solution satisfying the boundary conditions.  Moreover,
\eqref{eigreal} may have a solution for a finite number of negative
values of $b$ \cite[Section~3.2]{GS1}. So we still have an infinite
number of purely imaginary eigenvalues, but in addition to these
Kelvin waves there is also {\em continuous spectrum} filling the
interval where $0 \le b \le 1$.  Note that, if $0 < b < 1$, the
eigenvalue equation \eqref{eigreal} has a singularity at
$r = \bar r := \Omega^{-1}(b)$, which is referred to as a ``critical
layer'' in the physical literature. The interested reader is referred
to \cite{BG,FSJ,LL,RS} for a few recent contributions to the study of
Kelvin waves.

\subsection{Axisymmetric or Two-Dimensional Perturbations}
\label{sec23}

From now on we concentrate on the spectrum of the linearized operator
outside the imaginary axis. The stability equation \eqref{eigscalar}
is difficult to analyze in general, but important insight can be
obtained by considering some particular cases.

To begin with, we restrict our attention to {\em axisymmetric 
perturbations} for which $m = 0$. In that case, we have $\gamma(r) = s$ 
for all $r > 0$, so that \eqref{eigscalar} reduces to the simpler
equation
\begin{equation}\label{Axi}
  -\partial_r \partial_r^* u_r + k^2 \Bigl(1 + \frac{\Phi(r)}{s^2}
  \Bigr)u_r \,=\, 0\,, \qquad r > 0\,.
\end{equation}
The analogy with the Rayleigh-Taylor equation \eqref{RT6} is striking,
and we see that the Rayleigh function $\Phi$ in \eqref{Axi} plays the
exact role of the buoyancy term $N^2 = -g\bar\rho'/\bar\rho$ in
\eqref{RT6}. Following the same approach as in Section~\ref{sec11}, we
conclude that, if $\Phi$ is everywhere nonnegative, equation
\eqref{Axi} has no nontrivial solution satisfying the boundary
conditions when $\Re(s) \neq 0$. Moreover, if $\Phi(r) < 0$ for some
$r > 0$, Synge's argument \cite{GS1,Sy} shows that equation \eqref{Axi}
has a nontrivial solution for sequence of real eigenvalues $s_n \to 0$, 
so that the positivity of the Rayleigh function is a necessary and 
sufficient condition for stability in the axisymmetric case.

\begin{rem}\label{rem:energy}
The analogy between columnar vortices and shear flows in stratified 
fluids was already noticed by Rayleigh \cite{Ra2}, and can be 
roughly explained as follows. In a stratified fluid, exchanging 
the positions of two fluid particles located on the same vertical 
line results in a gain or a loss of potential energy, depending 
on whether the fluid density is decreasing or increasing upwards. 
The first situation is thus stable, and the second unstable. 
A similar effect occurs in vortices, even if the fluid is homogeneous, 
because the centrifugal force (which plays the role of gravity) 
varies as a function of the distance to the vortex center. It 
turns out that exchanging two fluid particles on the same radial
line results in a gain or a loss of energy depending on the 
sign of the Rayleigh function $\Phi$, and that a stable 
``stratification'' corresponds to $\Phi \ge 0$. 
\end{rem}

We next consider {\em two-dimensional perturbations}, which 
correspond to $k = 0$. In that case, the stability equation 
\eqref{eigscalar} reduces to 
 \begin{equation}\label{Twod}
  -\partial_r (r^2 \partial_r^* u_r) + m^2 u_r + \frac{imr W'(r)}{
  \gamma(r)}\,u_r \,=\, 0\,, \qquad r > 0\,.
\end{equation}
Here we can compare with the Rayleigh stability equation 
\eqref{SF3}, and we see that the vorticity derivative $W'$ 
in \eqref{Twod} plays the role of the second order derivative
$U''$ in \eqref{SF3}. Thus, proceeding as in Section~\ref{sec12}, 
we conclude that, if $W'$ does not change sign, equation 
\eqref{Twod} has no nontrivial solution satisfying the boundary 
condition if $\Re(s) \neq 0$. The monotonicity of the 
vorticity profile $W$ is thus a sufficient condition for 
stability with respect to two-dimensional perturbations, 
but as in the case of shear flows this condition is not 
necessary in general (and no sharp stability criterion is
known).

\begin{rems}\label{rem:part}\quad\\
{\bf 1.} For any localized vortex, the monotonicity of the 
vorticity distribution $W$ {\em implies} the positivity of 
the Rayleigh function $\Phi$. Indeed, if $W$ is monotone, 
then $W(r) \to 0$ as $r \to \infty$ (otherwise the vortex 
would not be localized), hence $W$ does not change sign, and 
the reconstruction formula
\begin{equation}\label{Omrep}  
  \Omega(r) \,=\, \frac{1}{r^2}\int_0^r W(s) s\dd s\,, \qquad r > 0\,,
\end{equation}
shows that $\Omega$ has the same sign as $W$. Thus $\Phi = 2\Omega W 
\ge 0$. \\[1mm]
{\bf 2.} In view of the previous remark, if we extrapolate the
conclusions obtained in the particular cases considered above, one
may be tempted to conjecture that a columnar vortex with monotone
vorticity distribution $W$ is (spectrally) stable for all values of
the Fourier parameters $m,k$. That daring claim has not been proved or
disproved so far, and it is good to keep in mind that, in the present
state of affairs, there is no analog of Squire's theorem for columnar
vortices.  In other words, there is no argument indicating that the
most unstable modes (if any) should always correspond to axisymmetric
or two-dimensional perturbations.
\end{rems}

\subsection{Howard Identities}\label{sec24}

We assume henceforth that $\Phi(r) > 0$ and $W'(r) < 0$ for all
$r > 0$, so that the vortex under consideration is stable with respect
to axisymmetric or two-dimensional perturbations. Our goal is now to
study the eigenvalue equation \eqref{eigscalar} in the general case
where $m \neq 0$ and $k \neq 0$. It is convenient to write the
spectral parameter as $s = m(a-ib)$, where $a,b \in \R$, so that
\begin{equation}\label{gamma1def}
  \gamma(r) \,=\, s + im\Omega(r) \,=\, im\ggamma(r)\,, 
  \qquad \hbox{where}\quad \ggamma(r) \,=\, \Omega(r) - b -ia\,.
\end{equation}
When $a \neq 0$, we have $\ggamma(r) \neq 0$ for all $r > 0$, and 
equation \eqref{eigscalar} can be written in the condensed form
\begin{equation}\label{eigscalar2}
   -\partial_r\bigl(\cA(r)\partial_r^* u_r\bigr) + \cB(r) u_r \,=\, 0\,,
   \qquad r > 0\,,
\end{equation}
where  $\partial_r^* = \partial_r + \frac1r$ and
\begin{equation}\label{ABdef}
  \cA(r) \,=\, \frac{r^2}{m^2 + k^2 r^2}\,, \qquad 
  \cB(r) \,=\, 1 + \frac{r}{\ggamma(r)}\partial_r\biggl(
  \frac{W(r)}{m^2+k^2 r^2}\biggr)
  - \frac{k^2}{m^2}\,\frac{\cA(r)\Phi(r)}{\ggamma(r)^2}\,.
\end{equation}

If we assume that \eqref{eigscalar2} has a nontrivial solution that is
regular at the origin and decays to zero at infinity, we can multiply
both sides of by $r\bar u_r$ and integrate over $\R_+$ to arrive at
the identity
\begin{equation}\label{HG0} 
  \int_0^\infty \Bigl(\cA(r)|\partial_r^* u_r|^2 + \cB(r) |u_r|^2 
  \Bigr)r \dd r \,=\, 0\,. 
\end{equation}
As $|\ggamma(r)| \ge |a| > 0$ for all $r > 0$, we deduce from \eqref{ABdef} that
\[
  |1 - \cB(r)| \,\le\, \frac{C}{m^2}\Bigl(\frac{1}{|a|} 
  + \frac{1}{|a|^2}\Bigr)\,, \qquad r > 0\,,
\]
for some constant $C > 0$ depending only on the vorticity profile
$W$. In particular, if we suppose that $|a| > M : = \max(1,2C)$, then
$\Re \cB(r) > 0$ for all $r > 0$, and taking the real part of
\eqref{HG0} we obtain a contradiction.  Thus equation
\eqref{eigscalar2} has no nontrivial solution if $|a| > M$.
Similarly, if we take the imaginary part of \eqref{HG0} use the
definitions \eqref{gamma1def}, \eqref{ABdef}, we obtain the relation
\begin{equation}\label{HG0im}
  a\int_0^\infty \biggl\{\frac{r}{a^2 + (\Omega{-}b)^2}\partial_r 
  \Bigl(\frac{W(r)}{m^2+k^2r^2}\Bigr) + \frac{2(b-\Omega(r))}{(a^2 
  + (\Omega{-}b)^2)^2}\,\frac{k^2}{m^2}\,\cA(r)\Phi(r)\biggr\}|u_r|^2 
  r\dd r \,=\, 0\,.
\end{equation}
If $a \neq 0$, the integral in \eqref{HG0im} must vanish. But the
first term in the integrand is negative since $W'(r)< 0$, and the
second one is negative too if we suppose that $b \le 0$, because
$\Omega(r) > 0$ for all $r > 0$. Thus we conclude from \eqref{HG0im}
that \eqref{eigscalar2} has no nontrivial solution if $a \neq 0$ and
$b \le 0$, see Fig.~1.

To obtain further information on the spectrum outside the imaginary
axis, we proceed as in the case of the Taylor-Goldstein equation
\eqref{TG3}, which was analyzed in Section~\ref{sec13}. Following
Howard's approach \cite{How,HG}, we first consider the differential
equation satisfied by the new function $w_r = u_r /
\ggamma(r)$.
Straightforward calculations that are reproduced in
\cite[Section~3.4]{GS1} show that $w_r$ satisfies
\begin{equation}\label{eqwr}
  - \partial_r \Bigl(\ggamma(r)^2 \cA(r) \partial_r^* w_r\Bigr) + 
  \cD(r) w_r \,=\, 0\,, \qquad r > 0\,,
\end{equation}
where
\[
  \cD(r) \,=\, \ggamma(r)^2  + 2r\ggamma(r)\partial_r\Bigl(
  \frac{\Omega(r)}{m^2+k^2r^2}\Bigr)- \frac{k^2}{m^2}\,\cA(r)\Phi(r)\,.  
\]
In particular, if we multiply \eqref{eqwr} by $r\bar w_r$, integrate the 
result over $\R_+$, and take the imaginary part, we obtain the relation
\begin{equation}\label{HG1}
  2a \int_0^\infty \biggl\{(b-\Omega(r))\Bigl(\cA(r)|\partial_r^* w_r|^2 +  
  |w_r|^2\Bigr)  - r\partial_r\Bigl(\frac{\Omega(r)}{m^2+k^2r^2}\Bigr)
  |w_r|^2 \biggr\}r\dd r \,=\, 0\,.
\end{equation}
The second term in the integrand is positive, because $\Omega'(r) < 0$, 
and the first one is positive too if we assume that $b \ge 1$, so that 
$b - \Omega(r) > 0$ for all $r > 0$. We thus conclude from \eqref{HG1} 
that equation \eqref{eqwr}, hence also equation \eqref{eigscalar2}, has 
no nontrivial solution satisfying the boundary conditions if $a \neq 0$ and 
$b \ge 1$, see Fig.~1.  

\figurewithtex 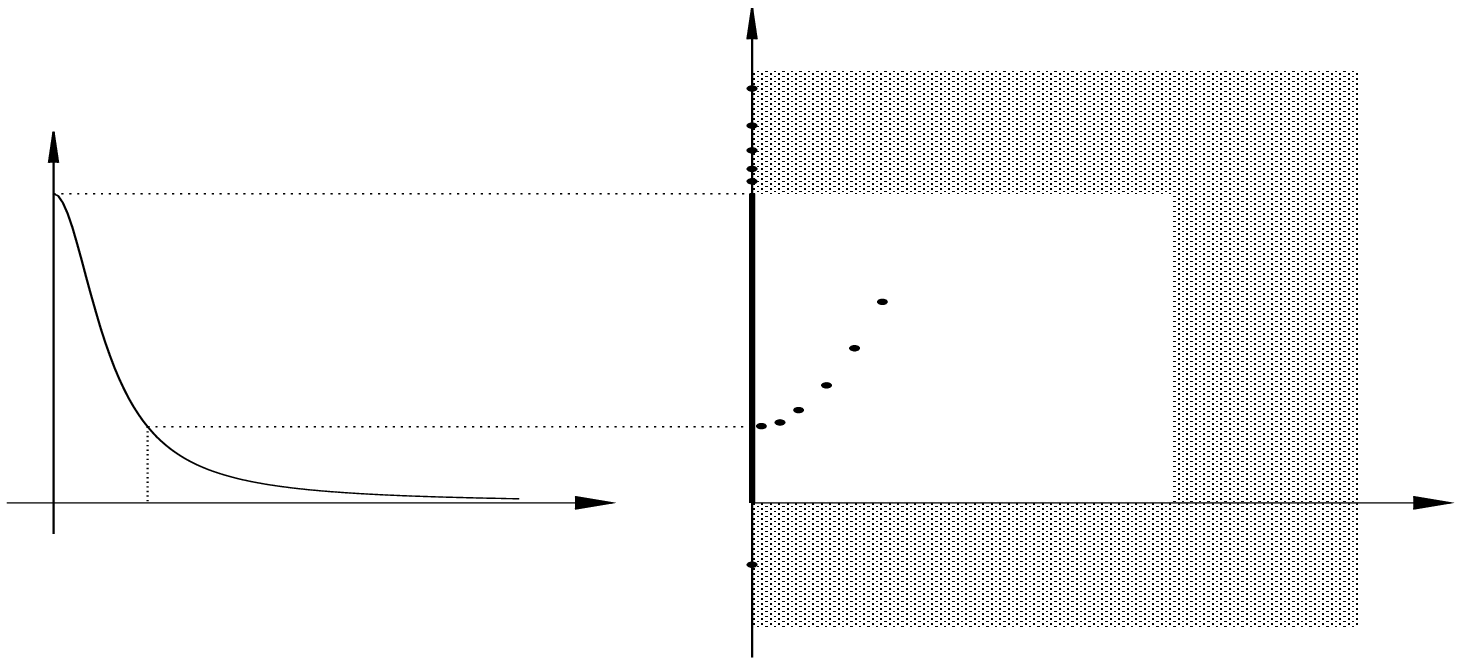 fig1.tex 6.50 15.00 {\bf Fig.~1:} 
The information obtained so far on the spectrum of the linearized
operator using the spectral parametrization $s=m(a-ib)$. Kelvin modes
are located on the imaginary axis $a = 0$, and accumulate only at the
upper edge of the essential spectrum, which fills the segment $a = 0$,
$b \in [0,1]$.  The rest of the spectrum, if any, consists of isolated
eigenvalues which are contained in the region $|a| \le M$,
$b \in [0,1]$, and can possibly accumulate only on the essential
spectrum. \cr

Next, we consider the function $v_r = u_r/\ggamma(r)^{1/2}$ which 
satisfies
\begin{equation}\label{eqvr}
  - \partial_r \Bigl(\ggamma(r) \cA(r) \partial_r^* v_r\Bigr) + 
   \cE(r) v_r \,=\, 0\,, \qquad r > 0\,,
\end{equation}
where
\[
  \cE(r) \,=\, \ggamma(r) + \frac{r}{2}\partial_r\Bigl(\frac{W(r) + 2\Omega(r)}
  {m^2+k^2r^2}\Bigr) + \frac14 \frac{\Omega'(r)^2}{\ggamma(r)} \,\cA(r)
  - \frac{k^2}{m^2}\,\frac{\cA(r) \Phi(r)}{\ggamma(r)} \,.
\]
If we multiply \eqref{eqvr} by $r\bar v_r$, integrate the result over $\R_+$, 
and take the imaginary part, we obtain the relation
\begin{equation}\label{HG1/2}
  -a \int_0^\infty \biggl\{\cA(r)|\partial_r^* v_r|^2 + |v_r|^2 +
  \frac{\cA(r)}{a^2 + (\Omega-b)^2}\Bigl(\frac{k^2 \Phi(r)}{m^2} 
  - \frac{\Omega'(r)^2}{4}\Bigr)|v_r|^2\biggr\}r\dd r \,=\, 0\,,
\end{equation}
which is analogous to identity \eqref{TG6}. Introducing the 
``Richardson number"
\begin{equation}\label{Ridef}
  \Ri(r) \,=\, \frac{k^2}{m^2}\,\frac{\Phi(r)}{\Omega'(r)^2}\,,
\end{equation}
we deduce from \eqref{HG1/2} that equation \eqref{eqvr}, hence also
equation \eqref{eigscalar2}, has no nontrivial solution satisfying the
boundary conditions if $a \neq 0$ and $\Ri(r) \ge 1/4$ for all
$r > 0$. Unfortunately, unlike for the Taylor-Goldstein equation, the
Richardson number \eqref{Ridef} depends on the Fourier parameters
$m,k$, and it is obvious that the inequality $\Ri(r) \ge 1/4$ cannot
hold for all values of $m$ and $k$. So the above approach fails to
give any stability criterion that would hold for arbitrary
perturbations. The situation is plainly summarized by Howard and 
Gupta in \cite{HG}\:

\begin{quote}
``The overall conclusion of this consideration of the non-axisymmetric
case is thus essentially negative\: the methods used to derive the
Richardson number and semicircle results in the axisymmetric case
reproduce the known results of Rayleigh for two-dimensional
perturbations and pure axial flow, but seem to give very little
more. In fact the present situation with regard to non-axisymmetric
perturbations seems to be very unsatisfactory from a theoretical point
of view."
\end{quote}

\begin{rem}\label{rem:semi}
In the spirit of Howard's semi-circle law for shear flows \cite{DH},
it is possible in the case of columnar vortices to locate the
(hypothetical) unstable modes in a slightly more precise way than what
is depicted in Fig.~1, see e.g. \cite{Eck}.  We do not comment further
on that, because in the next section we give conditions on the
vorticity profile which entirely preclude the existence of unstable
eigenvalues.
\end{rem}

\section{Spectral Stability of Inviscid Columnar Vortices}
\label{sec3}

In this section, we present the main results that were obtained recently 
in collaboration with D.~Smets \cite{GS1,GS2}. We first state our precise
assumptions on the unperturbed columnar vortex.

\smallskip\noindent{\bf Assumption H1:} {\em The vorticity profile 
$W : \overline{\R}_+ \to \R_+$ is a $\cC^2$ function satisfying $W'(0) = 0$, 
$W'(r) < 0$ for all $r > 0$, $r W'(r) \to 0$ as $r \to \infty$, 
and}
\begin{equation}\label{Winteg}
  \Gamma \,:=\, \int_0^\infty W(r) r\dd r \,<\, \infty\,.
\end{equation}

The crucial point here is the monotonicity of the vorticity
distribution $W$, which implies stability with respect to
two-dimensional perturbations, see Section~\ref{sec23}. We also
suppose that $W(r) \to 0$ as $r \to \infty$ fast enough so that the
integral in \eqref{Winteg} converges; in other words, the {\em total
  circulation} of the vortex is finite. It follows in particular that
$W(r) > 0$ for all $r > 0$, and the expression \eqref{Omrep} of the
angular velocity shows that $\Omega(r) > 0$ and $\Omega'(r) < 0$ for
all $r > 0$. As a consequence, the Rayleigh function $\Phi = 2\Omega W$ 
is positive everywhere, which implies stability with respect to 
axisymmetric perturbations too.

\smallskip\noindent{\bf Assumption H2:} {\em The ``Richardson function"
$J : \R_+ \to \R_+$ defined by
\begin{equation}\label{Jdef}
  J(r) \,=\, \frac{\Phi(r)}{\Omega'(r)^2}\,, \qquad  r > 0\,,
\end{equation}
satisfies $J'(r) < 0$ for all $r > 0$ and $rJ'(r)\to 0$ as $r\to \infty$.} 

\smallskip This second assumption is less natural, and probably only
technical in nature.  The quantity $J(r)$ appears in the definition of
the ``Richardson number" \eqref{Ridef}, which plays an important role
in the stability analysis of columnar vortices. If, for some given
value of the ratio $k^2/m^2 > 0$, the Richardson number \eqref{Ridef}
is not everywhere larger than $1/4$, assumption H2 implies the
existence of a unique $r_* > 0$ such that $\Ri(r) > 1/4$ if $r < r_*$
(stable region) and $\Ri(r) < 1/4$ if $r > r_*$ (possibly unstable
region). If we do not suppose that the function $J$ is monotone, more
regions have to be considered, which greatly complicates the
analysis. The monotonicity of $J$ is also essential to construct
simple subsolutions of equation~\eqref{eigreal} for large $r$, see
\cite[Section~4.6]{GS1}.  On the positive side, we emphasize that
assumptions H1 and H2 are satisfied in all classical examples, such as
the Lamb-Oseen vortex or the Kaufmann-Scully vortex. 

The following statement is our first main result. 

\begin{thm}\label{main1} {\bf \cite{GS1}}
Under assumptions H1, H2, the columnar vortex with vorticity profile 
$W$ is spectrally stable in the following sense. Given any $m \in \Z$ 
and any $k \in \R$ with $(m,k) \neq (0,0)$, the stability equation 
\eqref{eigscalar} has no nontrivial solution $u_r \in L^2(\R_+,r\dd r)$ 
if the spectral parameter $s$ has a nonzero real part. 
\end{thm}

Theorem~\ref{main1} asserts that, under assumptions H1, H2, the
linearized operator in \eqref{upert} has no unstable eigenmode of the
form \eqref{Four} with $\Re(s) \neq 0$ and $u \in L^2(\R_+,r\dd r)^3$.
In some sense, this answers a long-standing question dating back to
the pioneering contributions of Kelvin and Rayleigh. This rather
optimistic view has to be tempered for at least two reasons\: first,
the status of assumption H2 is unclear, and it is conceivable that the
conclusion of Theorem~\ref{main1} holds under the sole hypothesis that
the vorticity profile is monotone, although we do not know how to
prove that.  Next, the proof of Theorem~\ref{main1} given in \cite
{GS1} is very indirect, and does not give much insight into the
physical mechanisms leading to stability.  Therefore, it is not clear
if our approach can be applied to more complicated problems, such 
as the stability analysis of columnar vortices with nonzero
axial flow.

As is explained in Section~\ref{sec24}, if the angular Fourier mode
$m$ and the vertical wavenumber $k$ are both nonzero, the historical
approach to hydrodynamic stability based on integral identities such
as \eqref{HG0} does not seem sufficient to preclude the existence of
unstable eigenvalues in all regions of the complex plane, see Fig.~1.
However, it is easy to verify that all unstable eigenvalues (if any)
are simple, isolated, and depend continuously on the vortex profile
$W$, which can be considered as an infinite-dimensional parameter in
the differential equation \eqref{eigscalar}.  In addition, for the
rescaled Kaufmann-Scully vortex
\begin{equation}\label{Weps}
  W_\epsilon(r) \,=\, \frac{2}{(1+ \epsilon\,r^2)^2}\,, \qquad \hbox{where}
  \quad  0 \,<\, \epsilon \,\le\, \frac{4k^2}{m^2}\,,
\end{equation}
a direct calculation shows that the Richardson number \eqref{Ridef} satisfies 
$\Ri_\epsilon(r) \ge 1/4$ for all $r > 0$. By Howard and Gupta's result 
\cite{HG}, it follows that the associated linearized operator has no unstable 
eigenvalue in the Fourier subspace indexed by $m,k$. 

These observations suggest the following contradiction argument to
prove Theorem~\ref{main1}. Assume that, for some vorticity profile $W$
satisfying assumptions H1 and H2, the linearized operator in
\eqref{upert} has an unstable eigenmode of the form \eqref{Four} for
some $s \in \C \setminus i\R$ and some Fourier parameters $m \in \N$,
$k \in \R$. We know from the results of Section~\ref{sec23} that both
$m$ and $k$ are necessarily nonzero. The idea is now to perform a {\em
  continuous homotopy} $(W_t)_{t \in [0,1]}$ between the original
profile $W_0 := W$ and the reference profile $W_1 := W_\epsilon$,
where $W_\epsilon$ is defined in \eqref{Weps}.  For small $t$, the
linearized operator associated with $W_t$ has an unstable eigenvalue
$s(t)$ which depends continuously on $t$ and satisfies $s(0) = s$. But
we also know that, for $t = 1$, the linearized operator associated
with the reference profile $W_\epsilon$ has no unstable eigenvalue at
all.  Thus we logically conclude that there exists some
$t_* \in (0,1]$ such that the unstable eigenvalue $s(t)$ merges into
the continuous spectrum on the imaginary axis at $t = t_*$. The core
of our contradiction argument is the claim that, under assumptions H1
and H2, such a merger is actually impossible.

The way we actually arrive at a contradiction is not easily described in a 
few lines, and the interested reader is referred to \cite[Section~4]{GS1} 
for full details. If $t_n$ is an increasing sequence converging to
$t_*$, we denote $s_n = s(t_n) = m(a_n - i b_n)$, so that $a_n \to 0$
as $n \to \infty$ by construction. Also, extracting a subsequence if
needed, we can assume that $b_n \to \bar b \in [0,1]$ as
$n \to \infty$, see Fig.~1. For simplicity, we suppose here that
$0 < \bar b < 1$, but of course the limiting cases $\bar b = 0$ and
$\bar b = 1$ are also treated in \cite{GS1}.  If $u_r^n$ denotes the
(suitably normalized) eigenfunction associated with the eigenvalue
$s_n$ and the vorticity profile $W_{t_n}$, it is straightforward to
verify that $u_r^n$ converges as $n \to \infty$ to a solution $u_r$ of 
the limiting equation \eqref{eigreal}, where $b = \bar b$ and
$\Omega, W, \Phi$ denote the angular velocity, vorticity, and Rayleigh
function of the vortex profile at the bifurcation point $t =
t_*$. That equation has a singularity at the point
$\bar r := \Omega^{-1}(\bar b)$, and it is crucial to study the
behavior of $u_r$ in the vicinity of $\bar r$ (this is what is
referred to as a {\em critical layer analysis} in the physical
literature). If $\Ri(\bar r ) > 1/4$, it is relatively easy to obtain
a contradiction from identity \eqref{HG1/2}, because all main terms in
the integrand are positive in that case. If $\Ri(\bar r ) < 1/4$, a
contradiction can be obtained by a careful study of the solutions of
\eqref{eigreal} near the singularity, and by the construction of
appropriate subsolutions in the region where $r > \bar r$, 
see \cite{GS1}. 

\begin{rem}\label{rem:H2}
The argument we have just sketched requires that assumption H2 be
satisfied by the interpolated profile $W_t$ for all $t \in [0,t_*]$. 
For that reason, we cannot use a linear interpolation of the form
$W_t = (1-t)W + t W_\epsilon$, because the class of vorticity profiles
satisfying H2 is not a linear space nor even a convex set. Thus an
additional technical difficulty in our proof is the necessity of
constructing {\em ad hoc} interpolation and approximation schemes in
the nonlinear class of profiles satisfying assumption H2, see
\cite[Section~6.4]{GS1}.
\end{rem}

To state our second main result, we return to the linearized system 
\eqref{upert} which we write in condensed form $\partial_t \tilde u 
= L \tilde u$. The linearized operator $L$ is given by
\begin{equation}\label{Ldef}
  L \tilde u \,=\, \begin{pmatrix*}[l]
  -\Omega \partial_\theta \tilde u_r + 2 \Omega \tilde u_\theta 
  -\partial_r P[\tilde u] \\[1mm]
  -\Omega \partial_\theta \tilde u_\theta - W \tilde u_r  
  -\frac1r\partial_\theta P[\tilde u] \\[1mm]
  -\Omega \partial_\theta \tilde u_z -\partial_z P[\tilde u]\end{pmatrix*}\,,
\end{equation}
where $P[\tilde u]$ denotes the solution $\tilde p$ of elliptic 
equation \eqref{pressure}. Our goal is to solve the linearized system 
in the Hilbert space
\[
  X \,=\, \Bigl\{u = (u_r,u_\theta,u_z) \in L^2(\R^3)^3\,\Big|\, 
  \partial_r^* u_r + \frac1r \partial_\theta u_\theta + \partial_z u_z 
  = 0\Bigr\}\,,
\]
equipped with the standard $L^2$ norm.

\begin{thm}\label{main2} {\bf \cite{GS2}}
Assume that the vorticity profile $W$ satisfies assumptions H1, H2. 
Then the linear operator $L$ defined in \eqref{Ldef} is the generator of 
a strongly continuous group $(e^{tL})_{t \in \R}$ of bounded linear operators 
in the energy space $X$. Moreover, for any $\epsilon > 0$, there exists a 
constant $C_\epsilon \ge 1$ such that
\begin{equation}\label{eLbound}
  \|e^{tL}\|_{X \to X} \,\le\, C_\epsilon\,e^{\epsilon |t|}\,, \qquad
  \hbox{for all } t \in \R\,.
\end{equation}
\end{thm}

Estimate \eqref{eLbound} exactly means that the spectrum of the
evolution operator $e^{tL}$ in $X$ is contained in the unit circle of
the complex plane for all $t \in \R$. In that sense,
Theorem~\ref{main2} is arguably the strongest way of asserting that
the columnar vortex with vorticity profile $W$ is {\em spectrally
  stable}. In view of the Hille-Yosida theorem \cite{EN}, it follows
the spectrum of the generator $L$ is entirely contained in the
imaginary axis of the complex plane, and we have the following
resolvent bound for any $a > 0$\:
\begin{equation}\label{resbound}
  \sup\Bigl\{ \|(z - L)^{-1}\|_{X \to X} \,\Big|\, z \in \C\,,~|\Re(z)| 
  \ge a\Bigr\} \,<\, \infty\,.
\end{equation}
In fact, since $X$ is a Hilbert space, the Gearhart-Pr\"uss theorem 
\cite[Section~V.1]{EN} asserts that the resolvent bound \eqref{resbound} 
is {\em equivalent} to the group estimate \eqref{eLbound}. 

Let $L_{m,k}$ denote the restriction of the linearized operator $L$ to
the Fourier subspace indexed by the angular mode $m \in \Z$ and the
vertical wave number $k \in \R$. To prove Theorem~\ref{main2}, we fix
some spectral parameter $s \in \C$ with $\Re(s) = a \neq 0$ and we
consider the resolvent equation $(s - L_{m,k})u = f$, which is
equivalent to the system
\begin{equation}\label{eigsysf}
  \begin{array}{l}
  \gamma(r) u_r - 2\Omega(r)u_\theta \,=\, -\partial_r p + f_r\,, \\[1mm]
  \gamma(r) u_\theta  + W(r)u_r \,=\, -\frac{im}{r} p+f_\theta\,, \\[1mm]
  \gamma(r) u_z  \,=\, -ik p+f_z\,, \end{array} 
\end{equation}
where $\gamma(r) = s + im \Omega(r)$ and the pressure $p = P_{m,k}[u]$
is chosen so as to preserve the incompressibility condition
\eqref{incompODE}.  Our goal is to show that the solution of
\eqref{eigsysf} satisfies $\|u\| \le C(a)\|f\|$, where $C(a)$ is a
positive constant depending only on the spectral abscissa $a$; in
particular, the resolvent estimate is uniform in the Fourier
parameters $m,k$ and in the spectral parameter $s$ on the vertical
line $\Re(s) = a$. Such a uniform bound is essentially equivalent to
\eqref{resbound}, hence also to \eqref{eLbound} by the
Gearhart-Pr\"uss theorem.

If $(m,k) \neq (0,0)$, the resolvent system \eqref{eigsysf} can be
reduced to a scalar equation for the radial velocity $u_r$, which can
then be studied using the same techniques as in
Section~\ref{sec24}. This provides resolvent estimates with {\em
  explicit constant} $C(a)$ is some regions of the parameter space,
but that approach fails in other regions where we have to invoke a
contradiction argument that relies on the conclusion of
Theorem~\ref{main1}. Thus our proof is again non-constructive, and
does not provide any explicit expression for the constant $C(a)$ in
general. In particular, we do not know if $C(a) = \cO(|a|^{-N})$ as
$a \to 0$ for some $N \in \N$. Such an improved estimate would indicate
that the norm of the group $e^{tL}$ grows at most polynomially as
$|t| \to \infty$.

\section{Conclusion and Perspectives}
\label{sec4}

The results of the previous section apply to a large family of
columnar vortices, including all classical models in atmospheric flows
and engineering applications \cite{AKO,SS}. They provide the first
rigorous proof of spectral stability allowing for general
perturbations, without any particular symmetry. In this sense, they
solve an important problem that was formulated as early as 1880 by
Lord Kelvin in the pioneering work \cite{Ke}. However, many
interesting questions remain open\:

\begin{itemize}

\item Is assumption H2 really necessary for the conclusion of
  Theorem~\ref{main1} to hold\? Can one find a different proof, that
  does not rely on a non-constructive contradiction argument\?

\item Can one strengthen the conclusion of Theorem~\ref{main2} and
  show that the group norm $\|e^{tL}\|$ grows at most polynomially as
  $|t| \to \infty$\?

\item Is it possible to prove some spectral stability results for more
  general equilibria of the form $u = V(r) e_\theta + W(r) e_z$, which
  include a nonzero axial flow\?
 
\item Do our result give any useful information on the stability of
  columnar vortices in the slightly viscous case\?
 
\end{itemize}

\noindent
In a broader perspective, a long-term project is the stability
analysis of columnar vortices beyond the linear approximation, which
is a completely open problem in the absence of useful variational
characterization of such equilibria. In any case, we hope that our
contribution will serve as a starting point for new developments in
the stability analysis of concentrated vortices.

\end{document}